\documentclass[11pt]{article}

\usepackage{graphics}   
\usepackage{endnotes}            
\usepackage{amsfonts}            
\usepackage{amssymb}             
\usepackage{amsthm}              
\usepackage{amsmath}            
\usepackage{algorithm}           
\usepackage{algorithmic}
\usepackage[letterpaper]{geometry}
\usepackage{rotating}            

\usepackage{setspace}                
\usepackage{natbib}
\bibpunct{(}{)}{;}{a}{,}{;}

\begin{document}

\title{
A Dynamic Process Interpretation of the Sparse ERGM Reference Model\thanks{This work was supported by NSF award DMS-1361425.}
}

\author{
Carter T. Butts\thanks{Departments of Sociology, Statistics, and EECS and Institute for Mathematical Behavioral Sciences; University of California, Irvine; SSPA 2145; Irvine, CA 92697; \texttt{buttsc@uci.edu}}
}
\date{1/29/2018}
\maketitle

\begin{abstract}
Exponential family random graph models (ERGMs) can be understood in terms of a set of structural biases that act on an underlying reference distribution.  This distribution determines many aspects of the behavior and interpretation of the ERGM families incorporating it.  One important innovation in this area has been the development of an ERGM reference model that produces realistic behavior when generalized to sparse networks of varying size.  Here, we show that this model can be derived from a latent dynamic process in which tie formation takes place within small local settings between which individuals move.  This derivation provides one possible micro-process interpretation of the sparse ERGM reference model, and sheds light on the conditions under which constant mean degree scaling can emerge.\\[5pt]
\end{abstract}

\theoremstyle{plain}                        
\newtheorem{axiom}{Axiom}
\newtheorem{lemma}{Lemma}
\newtheorem{theorem}{Theorem}
\newtheorem{corollary}{Corollary}

\theoremstyle{definition}                 
\newtheorem{definition}{Definition}
\newtheorem{hypothesis}{Hypothesis}
\newtheorem{conjecture}{Conjecture}
\newtheorem{example}{Example}

\theoremstyle{remark}                    
\newtheorem{remark}{Remark}


\section{Introduction}

Stochastic models for network structure continue to advance rapidly both in terms of sophistication and empirical adequacy.  This advance provides important avenues for theoretical development, particularly in terms of identifying the formal relationships between classes of network models and the types of generative processes that can give rise to them.  In this paper, we consider a specific type of network model---the basic reference distributions that constitute the foundation on which exponential family random graph models are built---and the potential for such models to arise from a simple family of continuous time dynamic network processes.  In particular, we show how the reference family introduced by \citet{krivitsky.et.al:statm:2011} to allow generalization across vertex sets of varying size can be obtained as the limit of an unobserved dynamic processes in which nodes migrate among a set of social settings, forming ties only with those with whom they have a setting in common.

The remainder of the paper proceeds as follows.  After reviewing the notion of reference models in the ERGM context (and their interpretation) in this section, section~\ref{sec_scaling} introduces the substantively important problem of identifying network models with realistic mean degree scaling and the solutions that have been proposed in the literature.  As we show, the most obvious solution (degree constraints) leads to typically undesirable behavior, motivating the \citet{krivitsky.et.al:statm:2011} framework.  In section~\ref{sec_dyninterp}, we introduce a simple dynamic model whose equilibrium behavior approaches the Krivitsky model, providing an interpretation of the latter in terms of micro-level processes.  In addition to the asymptotic development, we employ simulation to probe the conditions that are required for the approximation to be effective in practice.  Section~\ref{sec_discuss} discusses additional issues related to time scales and provides an extension of the model to latent spatial processes, and section~\ref{sec_conclusion} concludes the paper.

\subsection{Background on ERGMs}

A random graph model in exponential family form (usually called an \emph{exponential family random graph model,} or ERGM) is a probability distribution for random graph $G$ on support $\mathcal{G}$ written as
\begin{equation}
\Pr(G=g|\theta,X) = \frac{\exp\left(\theta^T t(g,X)\right)}{\sum_{g'\in \mathcal{G}}\exp\left(\theta^T t(g',X)\right) h(g',X)} h(g,X) \label{e_ergm}
\end{equation}
where $X$ represents an arbitrary covariate set, $t:\mathcal{G},X \mapsto \mathbb{R}^k$ is a vector of \emph{sufficient statistics} governing the form of the distribution, and $\theta \in \mathbb{R}^k$ is a parameter vector.  The \emph{reference measure,} $h$, is a function taking $\mathcal{G},X$ into the (finite) non-negative real numbers; by definition, $h(g,X)=0$ for all $g\not\in\mathcal{G}$, and $h(g,X)>0$ for all $g\not\in\mathcal{G}$.  

Although ERGMs have many valid interpretations, it is perhaps most straightforward to think of them as parameterizing a random graph model in terms of a propensity to produce graphs with higher ($\theta_i>0$) or lower ($\theta_i<0$) values of the sufficient statistics, relative to some ``baseline'' reference model.  The nature of the reference model is determined by $h$, as can be appreciated by considering the form of Equation~\ref{e_ergm} as $\theta \to 0$.  In this limit, we obtain $\Pr(G=g|X) \propto h(g,X)$, demonstrating that $h$ is quite literally the distribution that arises (up to a constant of proportionality) when the social ``forces'' associated with $\theta$ are absent.  Although not emphasized in early work (indeed, many papers tacitly take $h(g,X)=1$ without mentioning it), interest in explicitly specifed reference measures as increased both due to the extension of ERGMs to valued settings \citep[where they are inescapable; see e.g.][]{krivitsky:ejs:2012,krivitsky.butts:sm:2017} and due to an increasing awareness that well-chosen reference measures are critical for both generalization of network models to new settings \citep{krivitsky.et.al:statm:2011,butts.almquist:jms:2015} and for ensuring good inferential properties \citep{schweinberger.et.al:tr:2017,krivitsky.kolaczyk:ss:2015}.

From a substantive standpoint, the reference measure can be viewed as expressing the effect of underlying (possibly unobserved) background conditions or processes on the structure of the graph.  The distinction between these factors and those associated with $t$ is somewhat arbitrary in the context of a model under fixed conditions (since it is always possible to fold elements of $h$ into $\theta^T t(G)$ by appropriate choice of $t$ and $\theta$),\footnote{This is conventionally done by the use of \emph{offset terms} \citep[see e.g.][]{krivitsky.kolaczyk:ss:2015} that deterministically alter elements of $\theta$.} but becomes important when generalizing to new conditions or when seeking to interpret $\theta$ (since the model parameters by construction measure deviations from what would be obtained from $h$ alone).  From this point of view, the reference measure encodes a form of prior knowledge regarding graph structure, whose choice should be motivated by the properties of the system being modeled.  This, in turn, requires an understanding of how particular choices of reference relate to underlying graph processes. 

\subsection{Reference Measures and their Interpretations}

Although it is tempting to regard the choice of reference measure as a purely technical device, its role in defining an underlying graph distribution on which other factors operate suggests a substantive interpretation.  More precisely, we say that a hypothetical social process supplies an \emph{interpretation} for a given ERGM family if realizations of networks arising from this process have a distribution corresponding to the ERGM in question.  It should be noted that determining whether a particular ERGM family admits a given interpretation depends on the type of ``realization'' involved; in different settings, one may be interested e.g. in cross-sectional observations of networks at random times, at aggregated tie structures from some known origin state over some interval, etc.  In our case, we focus on the context of a network-generating process whose state is observed at a random moment.  Numerous distinct processes may in some cases lead to a common distributional family, and hence a given ERGM family may have many interpretations - and, by turns, the mere fact that an ERGM family admits a particular process interpretation does not necessarily mean that empirical observations that are well-modeled by said family were generated by that process.  For instance, \citet{snijders:sm:2001} has shown that certain ERGM families can arise from a dynamic process of decision-making among boundedly rational agents.  This does not imply, however, that populations of networks that are well-modeled by those families must be the result of an underlying decision process.  Nevertheless, knowing that a particular ERGM family does (or does not) admit a decision-theoretic interpretation may provide useful clues regarding the social processes that could give rise to it; and, where there is \emph{a priori} reason to believe that a given network arose from a decision-theoretic process, such knowledge may serve as a powerful tool for model construction. 

While the question of interpretation can be applied to any ERGM family, it is particularly natural to seek interpretations of the reference models on which most ERGMs are built.  By a \emph{reference model}, we mean a choice of reference measure together with any statistics (i.e., elements of $t$) that are regarded as essential for any model within the family of interest.  For instance, the common (often implicit) choice of the counting measure ($h(g,X)\propto 1$) with no additional statistics leads to a reference family which is the uniform graph distribution on $\mathcal{G}$.  This class of reference families corresponds to the \emph{baseline models} of \citet{mayhew:jms:1984a}, which he defines to be the family of models having a uniform probability distribution on a specific state space (here, $\mathcal{G}$).  As Mayhew observes, the choice of state space is itself theoretical, and has substantial and non-trivial implications for social phenomena; he promotes the use of baseline models as substantively informed null hypotheses and as simple but possibly useful models in their own right.  Where $\mathcal{G}$ is the set of all graphs or digraphs on a given vertex set (as is typically the case), this reference model admits the interpretation of a draw from a stochastic process in which each edge arises independently and with probability 0.5.  In the ERGM setting, parameters for terms added to a uniform reference model can then be interpreted as forces biasing the distribution of graphs relative to the uniform graph distribution on $\mathcal{G}$.

Although the ``bare'' reference measure above leads to one class of reference models, one can define others.  For instance, in typical social network settings it is reasonable to assume that expected graph density will necessarily depart from the value of 0.5 implied by the uniform random graph, making this an unrealistic point of reference.  Instead, a more natural starting point is the general homogeneous Bernoulli graph, where ties arise independently with a fixed probability that is not necessarily equal to 0.5.  This reference model is obtained from the above by adding a single edge term, $t_e$, whose associated parameter $\theta_e$ is the logit of the expected density.  Parameters for terms added to the Bernoulli reference model can be the interpreted as forces biasing the graph distribution relative to a Bernoulli graph with parameter logit$^{-1}(\theta_e)$.  

Note that, since the uniform reference model is a special case of the Bernoulli reference model, the same ERGM family could be framed in terms of either reference (the key difference being in whether all other parameters are interpreted relative to the edge-only model, versus all parameters being interpreted relative to the zero-parameter model).  Nevertheless, the two reference models are distinct, and it is useful to draw a distinction between the two frames.  A pragmatic awareness of this distinction is also found in the standard advice to include an edge term in any ERGM fit to a network of non-trivial size (where a uniform model is not a plausible baseline): by always including an edge term (sometimes thought of as ``controlling for density''), one is ensuring that the general homogeneous Bernoulli graph remains the salient reference.

Obviously, many other reference models (and associated reference measures) are possible.  While several have been examined for valued \citep{krivitsky:ejs:2012} and rank-ordered \citep{krivitsky.butts:sm:2017} networks, relatively few have been investigated for the more common dichotomous case.  One of considerable importance is the ``sparse Bernoulli graph'' reference introduced by \citet{krivitsky.et.al:statm:2011} to create models that generalize across vertex sets of different sizes.  Since this model is central to our development, we consider it at greater length in the next section.

\section{Reference Models with Realistic Mean Degree Scaling} \label{sec_scaling}

It is a well-known observation that, in most social networks, changes in population size have at best a limited effect on mean degree.  For instance, human population has grown nearly 10-fold since the mid-18th century \citep[p42]{caselli.et.al:bk:2006}, but there is no evidence to suggest that the average number of sexual partners or close friends has shown a comparable increase.  Such observations lead to the common working assumption that, for a set of otherwise comparable networks, mean degree will be roughly constant in the size of the vertex set ($N$). 

Though seemingly unremarkable, this assumption is not in agreement with the common ERGM reference models described above; for both the uniform and Bernoulli reference models, mean degree must scale linearly in $N$.  This gives rise to the \emph{mean degree scaling problem}: generalizing from a model based on these references to a network of differing size from the one on which the model is based leads to mean degrees that do not scale correctly (at least, for typical social systems).  Note that this is a substantive theoretical issue, not solely a technical one.  The basic reference models fail because they (1) give every vertex an equal chance to be tied to every other vertex, and (2) hold that this chance does not depend on the size of the vertex set.  In such a world of unlimited tie formation opportunities, expected degree \emph{should} grow without bound; thus, the models are ill-behaved because they are based on assumptions that are typically unrealistic.  Correcting the problem requires modifying these assumptions.  There are, however, multiple ways to accomplish this, each of which has different consequences.  We briefly consider an important but problematic approach, before discussing the reference family on which we will focus in the balance of the paper.

\subsection{A State Space Solution: Constrained Maximum Degree}

Perhaps the most obvious explanation for asymptotically constant mean degree is the existence of bounds on the number of ties an individual can create or sustain.  Since the mean degree obviously cannot exceed the maximum degree, any such limits will necessarily constrain mean degree in the large-$N$ limit.  In our language, this alters the reference model by changing the state space, $\mathcal{G}$, from which $G$ is assumed to arise.  When the constant reference measure is retained and no other terms are required, the resulting reference model is a (Mayhewan) baseline model corresponding to a conditional uniform graph on the set of all graphs with maximum degree below the specified threshold value.  Adding an edge term leads to a Bernoulli-like reference model on the same support, though edges in this model are not entirely independent due to the constraint.

As noted, these reference families will exhibit constant mean degree scaling in the large-$N$ limit, and both implement the substantive mechanism of limited capacity to form and sustain ties.  However, this mechanism gives rise to a new problem: a reference model in which which edges form with fixed probability conditional on bounded maximum degree will \emph{saturate,} with almost all nodes eventually attaining the maximum degree.

To see how this effect arises, let us consider an arbitrary vertex from a reference model in which each pair of vertices is adjacent with probability $p$ conditional on both having less than $d_{\mathrm{max}}$ other edges, or else 0.  (The ERGM representation of such a model has a single edge term and uniform reference measure on the set of all graphs with maximum degree $d_{\mathrm{max}}$.)  Let $D_i$ be the degree of vertex $i$, and let $F_u$ be the fraction of vertices that are unsaturated (i.e., that have degrees less than $d_{\mathrm{max}}$).  For large $N$, the fraction of unsaturated vertices in $G\setminus i$ will be approximately the same as the entire graph.  Since each otherwise unsaturated vertex in the remainder of the graph is adjacent to $i$ with probability $p$ we have
 \begin{align*}
\Pr(D_i<d_{\mathrm{max}}) &\approx \sum_{i=0}^{\min\{F_u(N-1),d_{\mathrm{max}}-1\}} \binom{F_u(N-1)}{i} p^i (1-p)^{F_u(N-1)-i}\\
\intertext{Since the right-hand side is simply the cumulative distribution function of a Binomial$(F_u(N-1),p)$ random variable evaluated at $d_{\mathrm{max}}-1$, we may obtain an upper bound on it via Chernoff's inequality.  Noting that $F_u \approx \Pr(D_i<d_{\mathrm{max}})$ for randomly chosen $i$, we then have}
F_u &\le \exp\left(-\frac{1}{2p} \frac{\left(pF_u(N-1)-d_{\mathrm{max}}+1\right)^2}{F_u(N-1)}\right)\\
&= \exp\left(-\frac{p F_u(N-1)}{2} + \frac{d_{\mathrm{max}}-1}{2} - \frac{(d_{\mathrm{max}}-1)^2}{2 p F_u (N-1)}\right).\\
\intertext{Although we cannot solve for $F_u$ explicitly, the above is sufficient to show that no positive value of $F_u$ can be stable as $N \to \infty$.  To see this, let us hypothesize some fixed $F_u>0$, and consider the limit of the above as $N$ grows:}
&\xrightarrow[N\to \infty]{} \exp\left(-\frac{p F_u N}{2}\right),
\end{align*}
\noindent which converges to 0 for all $p$.  Thus, in a world where maximum degree is bounded and ties otherwise form independently with fixed probability, the fraction of unsaturated nodes must become vanishingly small as $N$ becomes large.

\begin{figure}[h]
  \centering
    \includegraphics[width=0.7\textwidth]{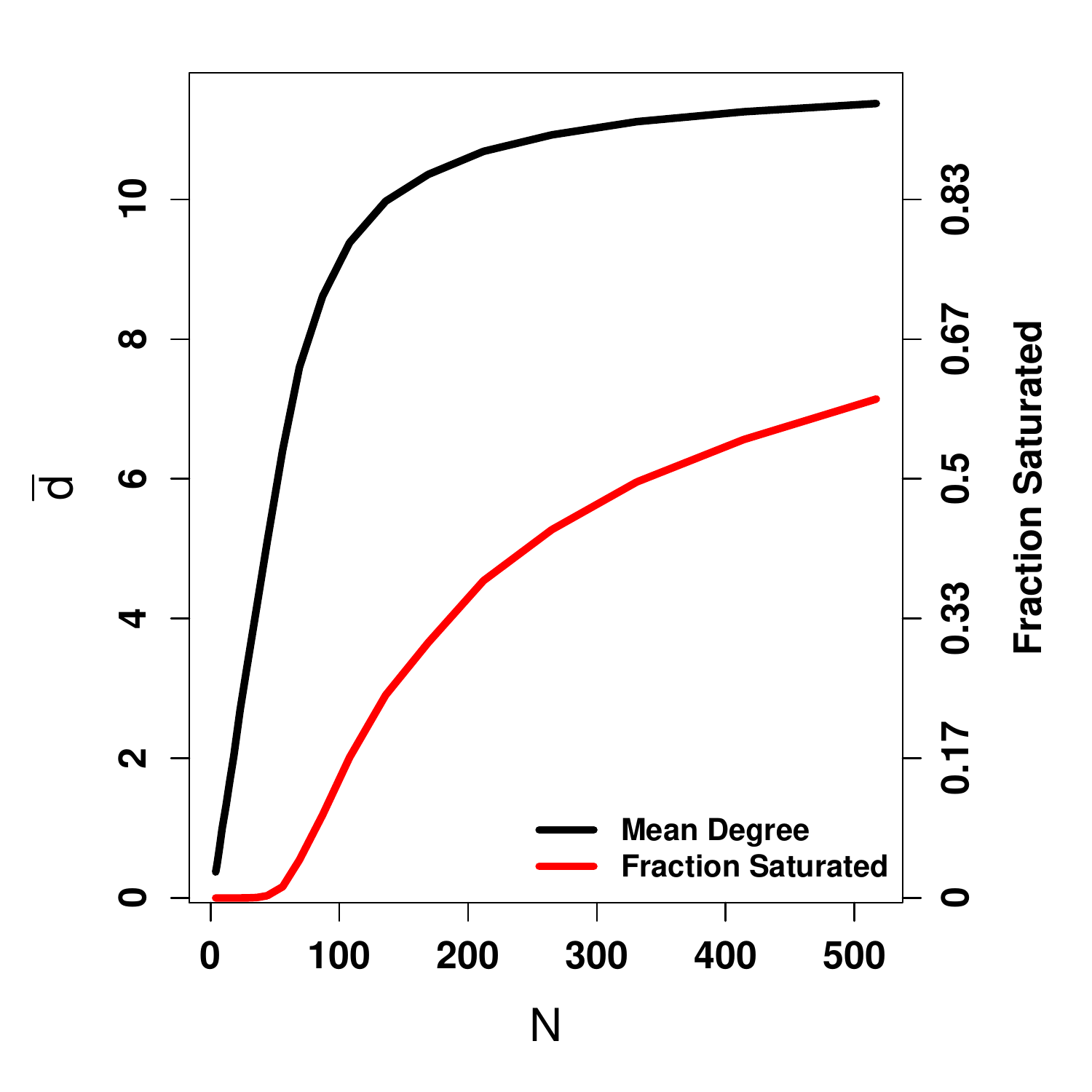}
    \caption{Degree saturation in a reference model with maximum degree 12, tie probability $p=0.12$.  As $N$ becomes large, mean degree (black line) approaches the maximum and an increasing fraction of nodes become ``saturated'' with edges (red line).\label{degsat}}
\end{figure}

This process is illustrated graphically in Figure~\ref{degsat}, which shows both mean degree and mean fraction of saturated nodes for a model with $p=0.12$ and $d_{\mathrm{max}}=12$.  (Simulations were performed using \texttt{ergm} \citep{hunter.et.al:jss:2008} for a constrained degree edge-only model with $\theta=\mathrm{logit}(0.12)$, with burnin values of $500 \tbinom{N}{2}$ and thinning intervals of $250 \tbinom{N}{2}$; 500 draws were obtained at each $N$.)  Mean degree grows linearly in $N$ before leveling off near $N=100$ (the point at which the expected degree would be equal to the degree constraint in an unconstrained Bernoulli graph), slowly approaching its limiting value of $d_{\mathrm{max}}$ as $N$ grows.  This leveling is accompanied by a corresponding rise in the fraction of saturated nodes, which begins to climb at around $N=50$ and grows steadly with $N$.  By $N=500$, around 60\% of all nodes are, on average, at maximum degree and unable to form additional ties.

This ``choking off'' of opportunities for tie formation is characteristic of large networks where degree scaling is held back by maximum degree constraints, and it may indeed be a real phenomenon for certain types of networks.  For most networks of interest, however, there is little evidence that the majority of nodes are saturated.  Individuals can and do form new friendships or make new professional contacts without being forced to jettison old ones, and the vast ranges in degree seen in acquaintenceship studies \citep{bernard.et.al:ch:1989,mccormick.et.al:jasa:2008} suggest that few if any persons are at the maximum capacity for such relationships to be sustained.  This suggests that the degree constraint model, while appealing in its simplicity, cannot be the explanation for mean degree scaling in most contexts.

\subsection{A Phenomenological Fix: the Krivitsky Reference Measure}

Another way to achieve constant mean-degree scaling is to return to the Bernoulli reference and simply force the expected density to fall as $N$ grows.  This approach, which was introduced by Krivitsky and colleagues (\citeyear{krivitsky.et.al:statm:2011}), can be summarized as follows.  Per \citet{butts:sm:2011b}, solving for the $\theta_e$ that leads to expected degree $\bar{d}\ge N-1$ gives us
\begin{align*}
\theta_e &= \log\left(\frac{\bar{d}}{N-1-\bar{d}}\right)\\
  &= \log \bar{d} - \log\left(N-1-\bar{d}\right).
\end{align*}
In the limit as $N/\bar{d} \to \infty$, we then have $\theta_e \to \log \bar{d} - \log N$.  This last was proposed by \citet{krivitsky.et.al:statm:2011} as the basis for an alternative reference model,
\begin{align}
\Pr(G=g|\bar{d}) &\propto \exp\left(\left(\log \bar{d} - \log N\right)t_e(g)\right) \nonumber\\
&= \exp\left(\theta_K t_e(g)\right) N^{-t_e(g)} \label{e_krivit}
\end{align}
where $\theta_K = \log \bar{d}$.  This model is obviously analogous to the Bernoulli graphs, but with reference measure $h(g,X) = N^{-t_e(g)}$.  We here refer to this as the \emph{Krivitsky reference measure}, and the associated model family as the sparse Bernoulli or Krivitsky reference model.

Motivated in this way, the Krivitsky reference model is a ``top-down'' solution to the mean degree scaling problem: having determined that mean degree should be asymptotically constant in $N$, we design a reference that has this property.  The resulting reference can be used to form ERGMs that generalize well across large differences in population size \citep{krivitsky.et.al:statm:2011}, have good inferential properties \citep{schweinberger.et.al:tr:2017}, and can be extended to accommodate arbitrary power-law degree scaling \citep{butts.almquist:jms:2015}.  It is reasonable, however, to ask how the edges within such a system ``know'' to become more improbable as the system size grows--i.e., to seek a microscopic interpretation of the reference model.  While such an interpretation is not necessary for the reference model to be gainfully employed, it may shed light on the emergence of constant degree scaling itself, and may further suggest conditions under which the Krivitsky reference model is likely to prove appropriate.  In the remainder of this paper, we will consider one such interpretation and investigate some of its implications for the relationship of underlying dynamic processes to cross-sectional network models.

\section{A Dynamic Interpretation of the Krivitsky Reference Model} \label{sec_dyninterp}

As noted, our core interest is in identifying a simple dynamic model whose equilibrium behavior is compatible with the Krivitsky reference model -- that is, a dynamic network model such that the distribution of a network observed at a random time will have a known ERGM form that matches an edge-only ERGM constructed using the Krivitsky reference measure (i.e., Equation~\ref{e_krivit}).  Since the process is intended to serve as a baseline relative to which other model terms are interpreted, it should be as simple and interpretable as possible, while allowing enough structure to generate the desired behavior (here, constant mean degree scaling).  Finally, the model should be plausible as a highly stylized social process.

The process we consider---which we will refer to simply as the ``contact formation model''--is motivated as follows.  As \citet{feld:ajs:1981} famously observed, formation of social ties often takes place within particular settings (which he dubbed ``foci'') in which (1) individuals encounter each other and (2) the social conditions needed for tie formation are present.  These conditions, like the foci themselves, may vary radically depending on the type of relation under consideration, and we do not attempt here to put any particular interpretation on them.  What we do assume, however, is that there exists a set of foci that can be considered to be exogenous and fixed on the time scale of network evolution, that the foci are distinct (i.e., each individual can occupy only one focus at a time), and that tie formation is only possible between individuals who occupy the same focus.  In particular, we assume that every pair of non-adjacent individuals occupying a common focus has some constant hazard, $r_f$ of forming a tie, while the hazard of tie formation for dyads whose endpoints reside in distinct foci is 0.  Once a tie is formed, however, foci are no longer needed for maintenance, and each currently existing edge is assumed to have a constant hazard $r_\ell$ of dissolution.  Finally, individuals are not assumed to occupy a given focus forever: we take each individual to have a constant hazard $r_m$ of migrating between foci.  When such a migration event occurs, the individual in question is assumed to choose a destination focus at random from the set of all foci (possibly returning immediately to the focus they were intending to depart).  

More formally, let us consider a set of $N$ vertices, $\mathcal{V}$, that are partitioned into a set of $M$ foci (each a set of vertices, some of which may be empty) $m_1,\ldots,m_M$.   Our dynamic model is a continuous time Markov process on state space $\mathcal{G}_\mathcal{V} \times \mathcal{F}_\mathcal{V}$, where $\mathcal{G}_\mathcal{V}$ is the set of all simple graphs on $\mathcal{V}$, and $\mathcal{F}_\mathcal{V} = \{1,\ldots,M\}^N$ is the set of all focus assignment vectors (i.e., vectors whose $i$th elements index the focus to which the $i$th vertex belongs).  Per the above, the possible transitions among elements of the state space consist of: (1) adding an edge between a non-adjacent pair sharing a common focus; (2) removing an edge; and (3) migrating a vertex from one focus to another.  The hazards for the events generating these transitions are as defined above (i.e., each pair at risk for tie formation forms a tie with instantanous rate $r_f$, each pair at risk for dissolution dissolves their tie with instantaneous rate $r_\ell$, and each individual has constant hazard $r_m/(M-1)$ of migrating from their current focus to each specific other focus).  This, and the piecewise constancy of the event hazards, makes it straightforward to calculate the state transition rates for the system as a whole (and thereby to simulate it).

Before deriving some properties of the behavior of this process, a few observations are in order.  First, we note that (trivially) the state space for this system is finite and has neither transient nor absorbing states.  Thus, it has a stationary distribution (or equilibrium), whose behavior will be of interest to us.  Second, the behavior of this process is quite dependent on the migration rate, $r_m$.  \emph{Our focus here is on the ``fast mixing'' regime, in which $r_m\gg r_f,r\ell$.}  That is, we will be interpreting our foci as intermittently occupied settings through which individuals move (e.g., locations at which they might encounter each other) rather than as long-term settings in which individuals are embedded.  Models based on slower-mixing regimes are possible, but do not lead to the same reference measure and are not our focus here.  Nevertheless, we do provide some results that suggest how these models would behave.

The model family used here is closely related to many other models that have been employed in the social network literature.  The constant hazard formation/dissolution framework arises as the continuous time limit of a separable temporal ERGM family \citep{krivitsky.handcock:jrssB:2014} with a single edge term for tie formation and dissolution (respectively).  This is the natural analog to the homogeneous Bernoulli graph for separable network dynamics, and its continuous limit is likewise a minimal model for separable network dynamics in continuous time.  Here, we elaborate that model by adding a simple form of latent block structure that affects only the tie formation process; models for latent block structure have a long history in network modeling \citep[see, e.g.][]{wang.wong:jasa:1987,snijders.nowicki:joc:1997,nowicki.snijders:jasa:2001}, though our interest is in mixing across realizations of block structure rather than inferring it.  (Additionally, the block structure is in our case ephemeral, whereas most prior work is focused on block structure that persists on the time scale of network evolution.)  Finally, we note that the continuous time structure employed here is closely related to the well-known stochastic actor-oriented models (SAOM) of \citet{snijders:jms:1996} and the longitudinal ERGMs of \citet{koskinen.et.al:ns:2015}.  Although we do not parameterize our process in the same manner as these families, it is possible to express the contact formation model in their respective languages (assuming that each is extended with the appropriate latent variable structure for foci, as is done e.g. in the SAB models of \citet{steglich.et.al:sm:2010}).  Our interest in the contact formation model is not hence in its novelty, but in its simplicity, close relationship to other well-known models, and potential to serve as a mechanistic underpinning for the use of the Krivitsky reference model.

\subsection{Asymptotic Development} \label{sec_asymp}

We now demonstrate that the Krivitsky reference model arises from random observations of the contact formation model in the large-$N$, high-$r_m$ limit.  Our development proceeds as follows.  First, we derive the mean degree of networks produced by the dynamic model in the regime of interest.  Next, we show that edges within observations from the dynamic model become asymptotically independent as $r_m \to \infty$.  Finally, we use these results to recapitulate the Krivitsky reference family, providing a direct interpretation of the edge parameter in terms of the parameters of the dynamic model.

\subsubsection{Limiting Mean Degree}

We derive the limiting mean degree by respectively deriving rates of edge formation and dissolution, and then solving for the resulting equilibrium.  To obtain the expected rate of edge formation, we first note that the expected number of dyads potentially at risk is equal to the number of vertex pairs sharing the same foci.  Since foci are exchangeable, this is in turn equal to the number of foci times the expected number of pairs within a given focus.  Since each vertex has an equal probability to be in each focus, the number within a given focus is binomially distributed with parameters $N$ and $1/M=P/N$; using the fact that the expected focus population is then $N/M=P$ and the expectation of the square of the focus population is likewise $N/M(N/M+1-1/M)=P(P+1-P/N)$, it follows that the expected number of same-focus pairs in the entire population is equal to
\[
\frac{N \left(P (P+1-P/N) -P \right)}{2 P} = \frac{P(N-1)}{2}.
\]
\noindent At any given time, not all of these pairs will be available for edge formation: some will already be adjacent.  Assuming a randomly mixed equilibrium state with mean degree $\bar{d}$ and density $\bar{d}/(N-1)$, the expected number of dyads at risk for tie formation will be approximately
\[
\frac{P(N-1)}{2} \left(1-\frac{\bar{d}}{N-1}\right)
\]
(this being exact in the limit as $r_m \to \infty$, where all focus memberships are instantaneously randomized).  The total rate of tie formation is then $r_f$ times the expected number of dyads at risk.

The situation with edge dissolution is more straightforward: in equilibrium, there are on average $\bar{d} N/2$ edges, each of which is lost at rate $r_\ell$.  In parallel to the above, the total rate of tie loss is equal to the loss rate times the number of edges at risk.

To solve for $\bar{d}$, we now note that, in equilibrium, the edge loss rate must equal the edge formation rate.  Thus
\begin{align}
r_\ell \frac{N}{2} \bar{d} &= r_f \frac{P(N-1)}{2} \left(1-\frac{\bar{d}}{N-1}\right) \label{e_gleq} \\
\bar{d} &= \frac{N-1}{N} \left[\frac{r_\ell}{r_f P}  + \frac{1}{N}\right]^{-1}  \label{e_mdegfull}\\
&\xrightarrow[N \to \infty]{} \frac{r_f P}{r_\ell}. \label{e_mdegasym}
\end{align}
Intuitively, this approximately the rate at which new edges would be gained by a single individual per edge lost, in the large-$P$ limit.  (Note that we have not assumed this in our development.)

It is interesting to compare the above with what would be seen in the low-migration limit.  As the migration rate approaches 0, the system can effectively be decomposed into a set of independent groups (one per focus) with purely internal formation/dissolution dynamics.  For an arbitrary $i$th group with size $P_i$, we can solve for the expected density $\delta_i$ by the equilibrium condition
\begin{align} 
r_f (1-\delta_i) P_i (P_i-1)/2 &= r_\ell \delta_i P_i (P_i-1)/2 \nonumber \\
r_f (1 - \delta_i) &= r_\ell \delta_i \nonumber\\
\delta_i &= \frac{r_f}{r_f+r_\ell}, \label{e_slowden}\\
\intertext{and hence the local equilibrium mean degree is}
\bar{d}_i &= (P_i-1) \frac{r_f}{r_f+r_\ell} \label{e_slowlocdeg}.
\end{align}
The fact that $\delta_i$ is constant can be leveraged to obtain $\bar{d}$ for the network as a whole.  By the exchangeability of foci, the expected number of edges for the entire graph (here denoted by $\mathbf{E} E$) is equal to $M \mathbf{E} E_i$, where $E_i$ is the edge count for an arbitrary focus.  Invoking the binomial distribution of focus sizes and the local density of Equation~\ref{e_slowden}, we have
\begin{align}
\mathbf{E} E &= M \mathbf{E} E_i \nonumber\\
&= M \sum_{j=0}^N \binom{n}{j} \left(\frac{1}{M}\right)^j \left(\frac{M-1}{M}\right)^{N-j} \frac{r_f}{r_f+r_\ell} \frac{j(j-1)}{2} \nonumber\\
&=\frac{r_f}{r_f+r_\ell} \frac{N(N-1)}{2M}. \nonumber
\end{align}
From this we obtain the global mean degree in the slow migration regime,
\begin{align}
\bar{d} &= \frac{2}{N} \mathbf{E} E \nonumber \\
&=\frac{r_f}{r_f+r_\ell} \frac{(N-1)}{M} \nonumber\\
&\xrightarrow[N \to \infty]{} \frac{r_f}{r_f+r_\ell} P \label{e_slowmdeg}
\end{align}
Observe that we recover the single-population mean degree for $M=1$, with $\bar{d} \to 0$ as $M \to \infty$.  (Note that the mean degree does remain positive for $M>N$, since no matter how many foci there are there is always the chance that some focus will have $P_i>1$. This becomes very rare for extremely large $M$, however, and the mean degree falls rapidly in this ``ultra sparse focus'' limit.)  Comparing the ``slow'' limit $\bar{d}$ of Equation~\ref{e_slowmdeg} with the ``fast'' limit $\bar{d}$ from Equation~\ref{e_mdegasym} gives us
\begin{align}
\frac{\bar{d}_{\mathrm{slow}}}{\bar{d}_{\mathrm{fast}}} &= \frac{r_f}{r_f+r_\ell} P \left[ \frac{r_f P}{r_\ell} \right]^{-1} \nonumber \\
&=  \frac{r_\ell}{r_f+r_\ell}, \label{e_slowfast}
\end{align}
\noindent from which a number of observations can be gleaned.  First, $\bar{d}_{\mathrm{slow}}<\bar{d}_{\mathrm{fast}}$, an effect that results from the crowding out of potential alters by existing adjacencies.  (Note that this effect would in principle vanish when the two models converge at $M=1$, but because of the assumptions involved in deriving $\bar{d}_{\mathrm{fast}}$ equation~\ref{e_slowfast} is only valid for $M>1$.)  Second, we also observe that the slow-fast gap is decreasing in $r_\ell$ and increasing in $r_f$.  This, too, is related to the ``choking'' effect.  To the extent that ties are ephemeral, the need to migrate to new local groups to find potential alters is diminished; relatedly, a high formation rate relative to $r_\ell$ results in a greater risk of ties ``piling up'' in one's local group and choking off opportunties to form more ties.  Where ties form quickly, last long, and can only be formed in small groups, migration is key to allow $\bar{d}$ to grow.

\subsubsection{Asymptotic Independence of Edges} \label{sec_indep}

In order to construct a reference model, we must know what terms should go in it; this amounts to knowing the conditional dependence structure of the edge variables in the observed graph \citep{pattison.robins:sm:2002}.  For the contact formation model, dependence among edges arises from a weak form of transitivity: if $i$ and $j$ are tied, then this implies they were likely in the same focus in the recent past; if $i$ is then also tied to $k$, this increases the conditional probability that $j$ and $k$ were likewise together, and hence able to form a tie.  This form of transitivity is strongest when $r_m \ll r_f,r_\ell$, as ties in this regime are essentially confined to vertices within common foci, and two-paths are very strong indicators of focus co-membership.  Intuitively, however, this effect is reduced as $r_m$ grows.  For $r_m \gg r_f,r_\ell$, vertex positions randomize on a time scale that is substantially faster than network evolution, giving co-membership little opportunity to affect structure.  Here we argue that, in the limit as $r_m \to \infty$ with $r_f/r_m, r_\ell/r_m \to 0$, this leads to asymptotic independence of edges.

We begin by considering the probability that two arbitrary vertices, $i$ and $j$ will be adjacent as a function of their past focus memberships.  For this purpose, it is convenient to reframe our network evolution process in terms of a series of \emph{formation events} and \emph{dissolution events} impacting a given vertex pair: dissolution events occur as a homogeneous Poisson process with rate $r_\ell$, while formation events occur as an inhomogeneous Poisson process with rate $r_f$ when $i$ and $j$ occupy the same focus and 0 otherwise.  Note that, in this representation, the events are not dependent on the state of the network -- formation events that occur when an edge is already present have no effect, and likewise with dissolution events that occur when no edge is present.  The probability of observing an $i,j$ edge at a random time is then the probability that the most recent event to occur was a formation event (since such an event would form an edge if one were not already present, and no dissolution event then occured before the system was observed).  Let $T_f$ be the time to the most recent formation event for $i$ and $j$, and $T_\ell$ the time to the most recent dissolution event.  Since we observe the system at a random time, we can exploit the memorylessness of the Poisson process to conclude that
\[
T_\ell \sim \mathrm{Exp}(r_\ell)
\]
(i.e., the time to the prior event is exponentially distributed with rate $r_\ell$.  It follows that the probability of observing an $i,j$ edge is
\[
\Pr(T_f<T_\ell) = \int_0^\infty r_\ell \exp(-r_\ell t) \Pr(T_f < t) dt.
\]
Unlike $T_\ell$, $T_f$ depends on focus co-membership; however, we can again invoke the memoryless property of the Poisson process to observe that $\Pr(T_f<t) = 1-\mathrm{Pois}\left(0\left|\int_0^t r_f S_x dx\right.\right)$ where $S_x =1$ if $i$ and $j$ occupy the same focus $x$ time units before the observation point and 0 otherwise.  Substituting for the Poisson pmf then gives us
\begin{align}
\Pr(T_f<T_\ell) &= \int_0^\infty r_\ell \exp(-r_\ell t) \left(1-\exp\left(-\int_0^t r_f S_x dx\right)\right) dt \nonumber\\
&=1-\exp\left(-r_\ell t\right) - \int_0^\infty r_\ell \exp\left(-r_\ell t -r_f C_t\right) dt \label{e_pegct}
\end{align}
where $C_t = \int_0^t S_x dx$ is the accumulated time during the prior $t$ time units in which $i$ and $j$ were occupying the same focus.  Although we cannot obtain $\Pr(T_f<T_\ell)$ without knowing $C_t$, we may immediately observe that since $0 \le C_t \le t$,
\[
0 \le \Pr(T_f<T_\ell) \le 1-\frac{r_\ell}{r_f+r_\ell}
\]
and, importantly, $c'_t \to c_t$ implies that $\Pr(T_f<T_\ell|C_t=c_c') \to \Pr(T_f<T_\ell|C_t=c_t)$.  In particular, this implies that if $C_t$ is constant across all possible states of the non-$i,j$ edge variables then the $i,j$ edge probability is unaffected by the graph state; where this holds for all arbitrary $i,j$, it follows that all edges will be independent.

Now let us consider $C_t$.  For an arbitrary vertex pair, we may without loss of generality fix one of them as a reference, and then construct $S_t$ in terms of periods in which the other party is or is not in the current reference focus.  Migration events in this representation occur as a Poisson process with rate $2r_m$, the rate doubling arising from our choice of representation (i.e., if the reference vertex moves, the reference location moves with him or her).  Since each event has probability $1/M$ of placing the vertices within the same focus (irrespective of starting position), such events themselves occur as a Poisson process with rate $2 r_m/M$.  During some fixed period time of length $t$, the number of events placing our vertex pair together is thus distributed as Pois$(2t r_m/M)$.  Each such event is followed by a period of co-residence, truncated by the next migration event (recall that events bringing a vertex back to the same focus are counted here as leading to separate, ``back-to-back'' intervals).  It follows from the properties of the Poisson process that such intervals have iid length distributions Exp$(2r_m)$.  %
%
%
Since the sum of iid exponential distributions is gamma distributed, it follows that $C_t$ is a Poisson$(2 t r_m/M)$ mixture of gamma deviates, the $i$th of which has distribution Gamma$(i,2 r_m)$ with expectation $i/(2 r_m)$ and variance $i/(2 r_m)^2$.  From this we can obtain the moments
\begin{align}
\mathbf{E} C_t &= \sum_{i=0}^\infty \mathrm{Pois}(i|2 t r_m/M) \frac{i}{2 r_m} \nonumber \\
&= \frac{t}{M} \nonumber
\end{align} 
and
\begin{align}
\mathrm{Var}(C_t) &= \sum_{i=0}^\infty \mathrm{Pois}(i|2 t r_m/M) \left[\left(\frac{i}{2 r_m}-\frac{t}{M}\right)^2+\frac{i}{(2 r_m)^2}\right] \nonumber \\
&= \sum_{i=0}^\infty \frac{(2 t r_m/M)^i \exp(-2 t r_m/M)}{i!} \frac{i + \left(i-2 t r_m/M\right)^2}{4 r_m^2}\nonumber\\
&=\frac{t}{r_m M} \label{e_varct}
\end{align} 
using standard properties of discrete mixtures.  We note that the expectation of $C_t$ does not depend on $r_m$; intuitively, it arises from an even division of the total available time ($t$) among all $M$ foci that the non-reference vertex could occupy.  The variance of $C_t$, however, \emph{does} depend on $r_m$, and in particular Var$(C_t) \to 0$ as $r_m \to \infty$.  An immediate consequence is that $C_t$ approaches a degenerate distribution at its expectation, taking the value $t/M$ almost surely as $r_m \to \infty$.

We are now in a position to establish the independence of edges under the contact formation model in the high-$r_m$ limit.  Let $y$ be the adjacency matrix for a realization of the graph structure $G$ under the contact formation model at a random time (with $Y$ the random variable equivalent), and let $i,j$ be an arbitrary edge pair. Denote by $Y^c_{ij}$ the set of all edge variables other than $\{i,j\}$.  Clearly, if $\Pr(Y_{ij}=1|Y^c_{ij}=y^c_{ij})=\Pr(Y_{ij}=1|Y^c_{ij}=(y')^c_{ij})$ for arbitrary alternative graph realizations $y'$, then $Y_{ij}$ is independent of the rest of the graph; if this is true for all $i,j$ pairs, then all edges are independent.  From Equation~\ref{e_pegct}, it follows that $\Pr(Y_{ij}=1|Y^c_{ij}=y^c_{ij}) \neq \Pr(Y_{ij}=1|Y^c_{ij}=(y')^c_{ij})$ only if $C_t | Y^c_{ij}=y^c_{ij} \neq C_t | Y^c_{ij}=(y')^c_{ij}$.  From Bayes's Theorem,
\[
p(C_t | Y^c_{ij}=y^c_{ij}) \propto \Pr(Y^c_{ij}=y^c_{ij} | C_t) p(C_t),
\]
where $p(C_t)$ is the marginal distribution of $C_t$.  It follows from Equation~\ref{e_varct}, however, that $p(C_t) \to 0$ for all $C_t \neq t/M$ as $r_m \to \infty$, and hence (invoking the fact that $\Pr(Y^c_{ij}=y^c_{ij} | C_t=t/M)$ is non-vanishing in $r_m$) $C_t | Y^c_{ij}=y^c_{ij}$ approaches a degenerate distribution at $t/M$ almost surely.  Since $t/M$ is the same for all values of $y^c_{ij}$, and since this applies to all $i,j$ pairs, it follows that all edge variables become independent in the contact formation model in the limit of fast migration.

\subsubsection{Constructing the Reference Model}

We now have the results that are required to construct a reference model for a random-time observation of a network from the contact formation model in the fast migration limit.  From Section~\ref{sec_indep}, we know that our model must be a homogeneous Bernoulli graph. Using the identity that the density is equal to $\bar{d}/(N-1)$, Equation~\ref{e_mdegasym} gives us the limiting graph density under the dynamic model.  To create an ERGM reference model that preserves this property, we employ a Bernoulli graph whose expected density matches the target.  In ERGM form, such a model has one statistic (the edge count, $t_e$), with a corresponding parameter $\psi$ equal to the logit of the density.  From Equation~\ref{e_mdegasym} we have
\begin{align}
\psi &= \log \frac{\frac{r_f P}{r_\ell (N-1)} }{ 1 - \frac{r_f P}{r_\ell (N-1)} } \nonumber\\
     &= -\log \left[ \frac{r_\ell (N-1) }{r_f P} - 1\right]. \nonumber\\
     &\approx -\log \left[ \frac{r_\ell N }{r_f P} -1\right] \nonumber\\
\intertext{for large $N$.  In large/sparse regime, $r_\ell N  \gg r_f P$, and thus}
     &\approx \log \frac{r_f P}{r_\ell} - \log N. \label{e_psi}
\end{align} 

$\psi$ thus has two components, one involving the relative rate of local edge formation versus deletion and the second a constant term involving $N$.  Writing the resulting reference model in ERGM form we have
\begin{align}
\Pr(G=g|\psi) &= \frac{\exp\left(\psi t_e(g)\right)}{\sum_{g' \in \mathcal{G}}\exp\left(\psi t_e(g')\right)} \nonumber\\
&= \frac{\exp\left[\left(\log \frac{r_f P}{r_\ell} - \log N\right) t_e(g)\right]}{\sum_{g' \in \mathcal{G}}\exp\left[\left(\log \frac{r_f P}{r_\ell} - \log N\right) t_e(g')\right]} \nonumber\\
&= \frac{\exp\left[\left(\log \frac{r_f P}{r_\ell}\right) t_e(g)\right]}{\sum_{g' \in \mathcal{G}}\exp\left[\left(\log \frac{r_f P}{r_\ell}\right) t_e(g')\right] N^{-t_e(g')}} N^{-t_e(g)} \nonumber\\
&= \frac{\exp\left[\theta_e t_e(g)\right]}{\sum_{g' \in \mathcal{G}}\exp\left[\theta_e t_e(g')\right] h_e(g')} h_e(g) \label{e_ergmbase} 
\end{align}
where $\theta_e = \log \tfrac{r_f P}{r_\ell}$ is the usual edge parameter, and $h_e(g)=N^{-t_e(g)}$ is the familiar Krivitsky reference measure.  Thus we are able to recapitulate the Krivitsky reference from the contact formation model, demonstrating that the latter can provide one interpretation for the former.

\subsection{Comparison with Simulated Dynamics}

Section~\ref{sec_asymp} provides an asymptotic development for an ERGM reference model from the contact formation process, but assumes that (1) the latent migration process is time-scale separated from the network dynamics, (2) the graph is sparse, and (3) $N$ is sufficiently large.  How sensitive is the reference model of equation~\ref{e_ergmbase} to these assumptions?  To investigate this, we compare the asymptotic mean degrees of equations~\ref{e_mdegfull} and \ref{e_mdegasym} to the mean degrees obtained by explicit simulation of the latent dynamic process.  We also compare the extent of triangulation seen under the contact formation model with what would be seen from the Bernoulli graph obtained in the fast-mixing limit.  To the extent that the process behavior mirrors the analytical form on which equation~\ref{e_ergmbase} is based, the corresponding reference model will serve as a good approximation for what would be observed under the contact formation process.

\subsubsection{Simulation Design}

To evaluate the behavior of the dynamic model, explicit continuous time simulations were performed for networks with varying values of $N$, $M$, and $r_m$.  To fix the (otherwise arbitrary) time scale, we chose the convention $r_f=1$; thus, all rates can be thought of as being expressed relative to the rate of tie formation.  To preserve constant mean degree scaling, we parameterized our simulations in terms of the ``expected local population'' $P=N/M$.  We also held the relative formation and dissolution rates constant, choosing $r_\ell=5$ to model a typical sparse graph regime.  For conditions in which $M=N/P$ was not an integer, model replications were randomly assigned to $M=\lceil N/P \rceil$ or $M=\lfloor N/P \rfloor$ such that the average value of $P$ was preserved across replications.  Each replication consisted of an independent simulation of 25 time units in length, seeded with an independent Bernoulli graph with the asymptotic mean degree, from which the final graph state was kept as the simulation outcome.  The 25 unit time period was chosen to be long enough to ensure convergence to equilibrium behavior, as assessed via a number of pilot runs.

Model parameters were varied as a full factorial design, with $N\in (50,100,200,400,800,1600)$, $P \in (5,10,20,40)$, and $r_m \in (1, 10, 25, 50, 100)$.  For each condition, 1,000 independent replications were performed, for a total of 120,000 simulation runs.  Simulation was performed using a custom script in the \textsf{R} statistical computing platform \citep{rteam:sw:2018}, using tools from the \texttt{statnet} library \citep{handcock.et.al:jss:2008,butts:jss:2008b}.

To obtain additional resolution regarding the transition from the slow mixing regime to the fast mixing regime, an additional simulation experiment using the above protocol was performed for $N=500$ with $P \in (5,10,20)$ and $r_m \in (1/125,1/25,1/5,1,5,25,125)$.  As above, 1,000 independent replications were performed per condition, for a total of 21,000 simulation runs.

\subsubsection{Recapitulating Mean Degree}

\begin{figure}[h]
  \centering
    \includegraphics[width=0.9\textwidth]{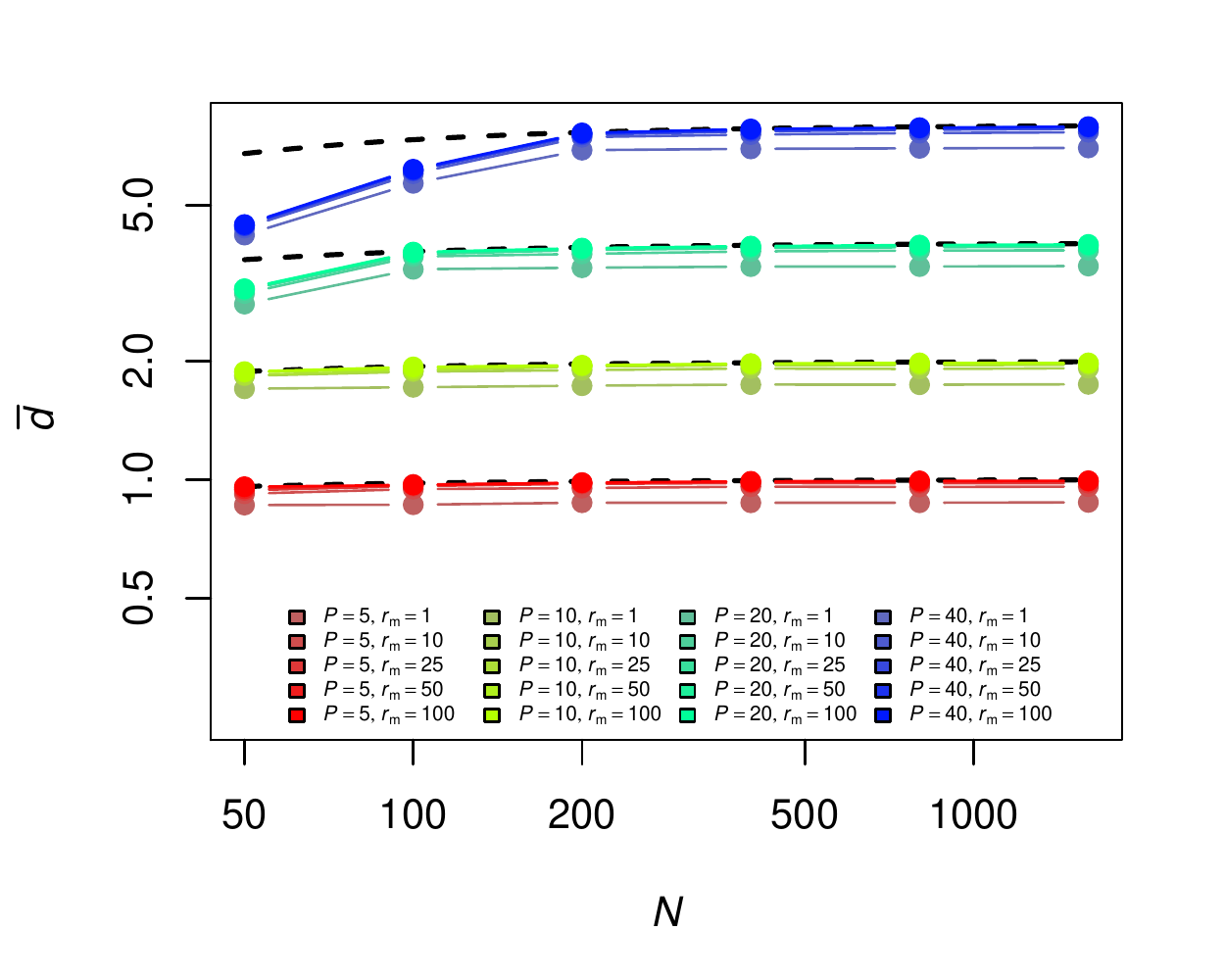}
    \caption{Expected degree, by $N$, $P$, and $r_m$.  For all local population sizes, mean degrees converge to their reference values (black dotted lines) as $N$ and $r_m$ become large.  (Note: 95\% confidence intervals are smaller than plotting symbols.) \label{f_convergence}}
\end{figure}

Our core question is whether the model family associated with the Krivitsky reference measure can accurately recapitulate the real behavior of the contact formation model.  As Figure~\ref{f_convergence} shows, the answer is affirmative: when movement rates are high, mean degrees under the dynamic model readily converge to the asymptotic limit for reasonable choices of $N$.  The network size required for accurate approximation increases with the asymptotic mean degree, reflecting the assumption of sparsity made in our development; in practice, $N\ge 25 \bar{d}$ seems to be sufficient for this condition to hold.

The more critical condition for convergence to the reference model is that the underlying migration rate, $r_m$, be high relative to the rate of edge turover.  As Figure~\ref{f_convergence} shows, movement rates on the same time scale as formation ($r_m=1$) leads to mean degree values that are systematically lower than predicted for all $N$ (as expected based on Equation~\ref{e_slowfast}). When migration rates are low, vertices spend considerable time ``trapped'' in foci with others to whom they are already tied, reducing their opportunities to form new ties; as migration rates grow relative to $r_f$, the number of local alters at risk for tie formation grows until it reaches its asymptotic limit of $(P-1)(1-\bar{d}/(N-1))\approx (P-1)(1-\bar{d}/N)$.  

\begin{figure}[h]
  \centering
    \includegraphics[width=0.9\textwidth]{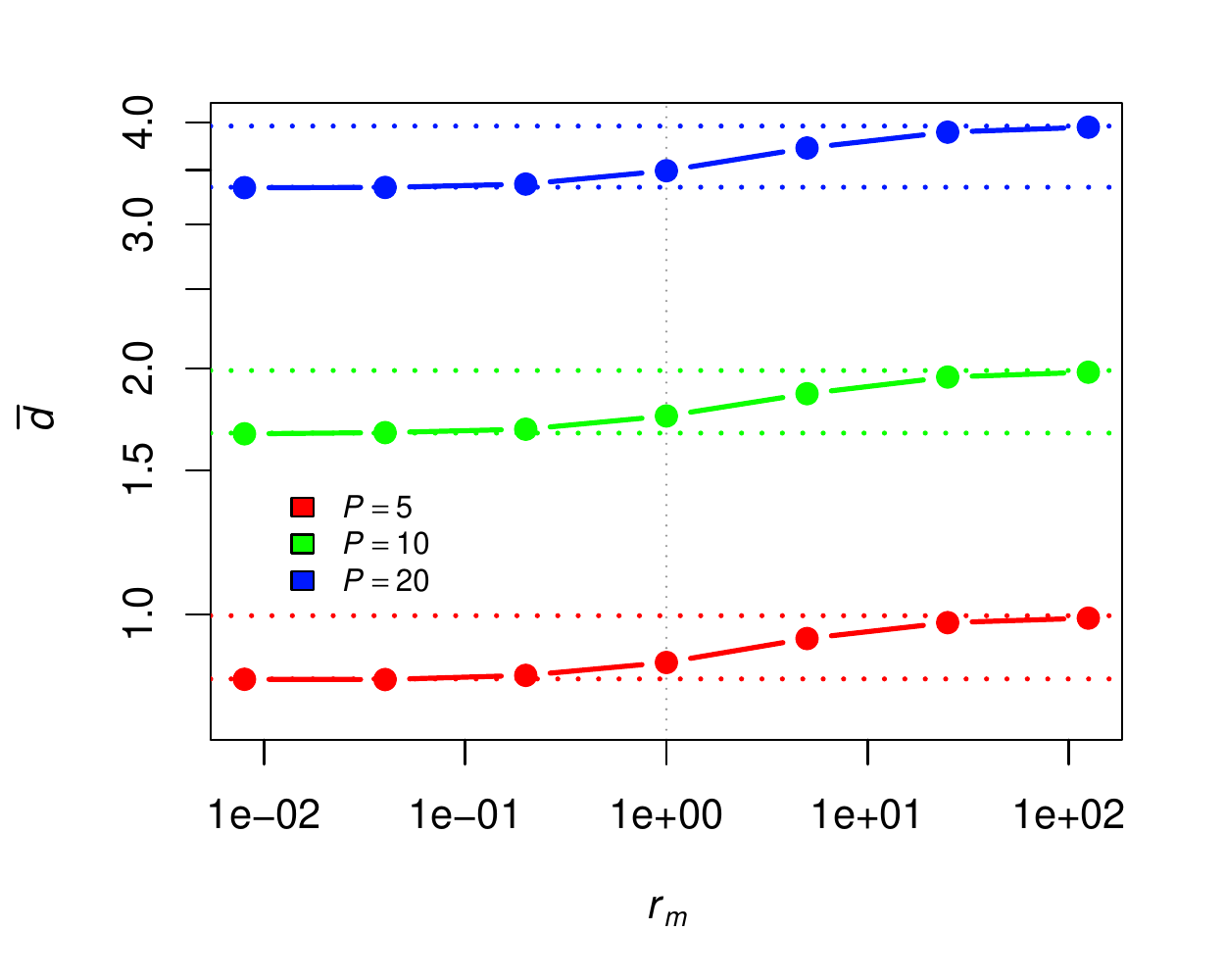}
    \caption{Expected degree at $N=500$, by $P$, and $r_m$.  Transition from slow mixing limit (lower dotted lines) to fast mixing limit (upper dotted lines) follows a characteristic form that is invariant to $P$.  (Note: 95\% confidence intervals are smaller than plotting symbols.) \label{f_migration}}
\end{figure}

This is shown in more detail in Figure~\ref{f_migration}, which shows the transition in $\bar{d}$ from $\bar{d}_{\mathrm{slow}}$ (lower dotted lines) to $\bar{d}_{\mathrm{fast}}$ (upper dotted lines) for $N=500$. This transition has a characteristic flattened sigmoidal form that does not depend on $P$, with effective time scale separation occuring when $r_m$ is between one and two orders of magnitude faster or slower than $r_f$.  While the quality of the fast-regime asymptotic approximation continues to improve as $r_m$ increases, our simulations suggest that $r_m$ should be at least an order of magnitude faster than $r_f$ (and ideally as much as two orders of magnitude faster) for the reference model to serve as a close approximation.

\subsubsection{Triangulation Rates}

\begin{figure}[h]
  \centering
    \includegraphics[width=0.9\textwidth]{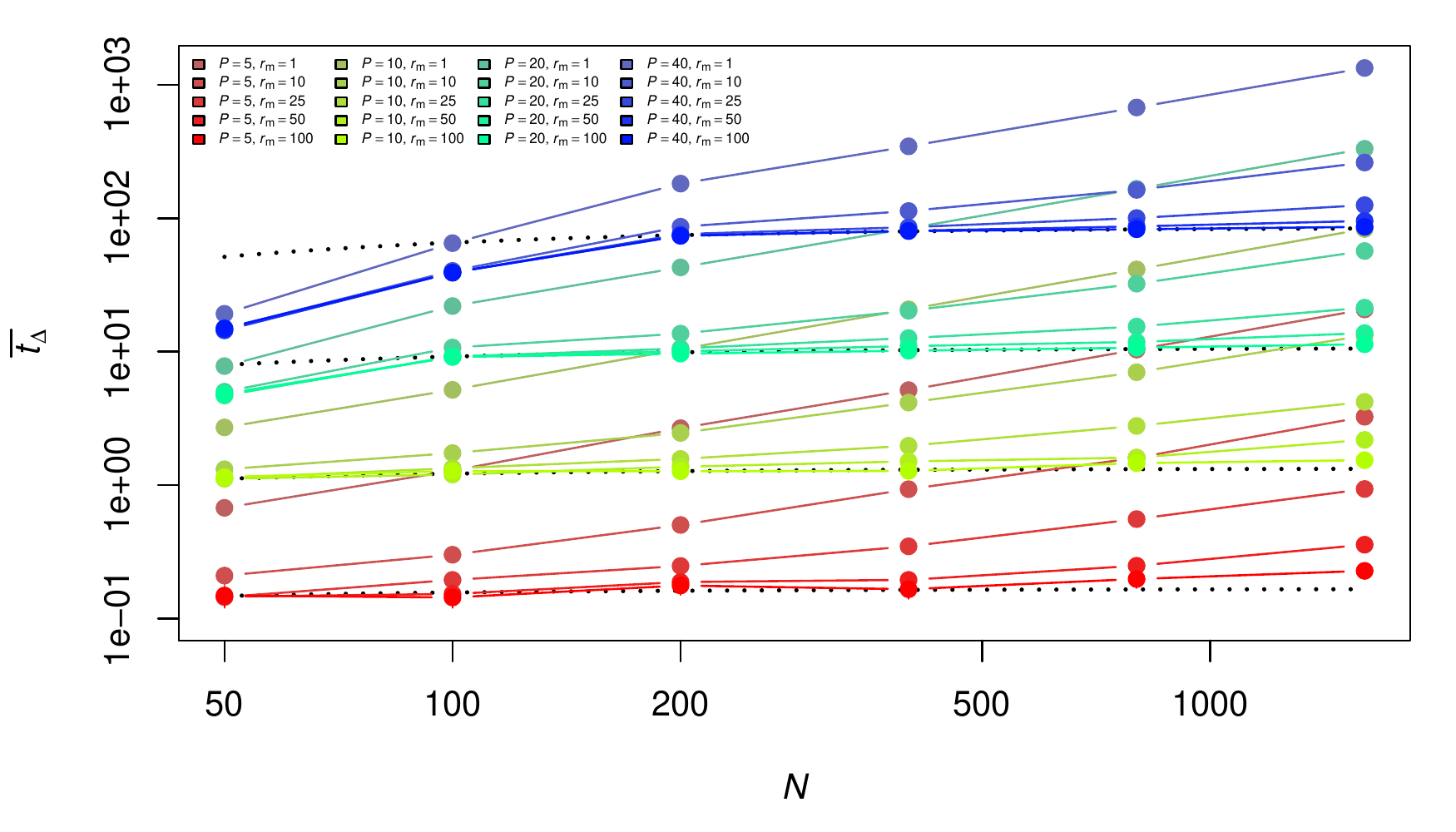}
    \caption{Expected triangle count, by $N$, $P$, and $r_m$.  For all $P$, mean values of $t_\Delta$ converge to their predicted values (black dotted lines) for sufficiently large $r_m$; however, higher values of $r_m$ are required for large $N$ and small $P$ to obtain convergence.\label{f_triangles}}
\end{figure}

Although we have proved that edges are independent in the high-migration limit, the trapping of nodes within foci at lower migration rates will lead to higher levels of triangulation than would be observed from an equivalent Bernoulli process.  How sensitive is the Bernoulli character of the marginal network distribution to the migration rate?  Figure~\ref{f_triangles} shows the mean triangle count in the simulated networks, versus their expectation under the reference model.  While low values of $r_m$ lead to very high levels of excess triangulation (as would be expected from a slow migration process), the behavior of the dynamic model converges to the Krivitsky reference as $r_m$ increases.  In particular, we find very close agreement over a broad range of network sizes when migration rates are 1.5-2 orders of magnitude faster than the tie formation rates.  This result is reassuring, in that it suggests that the Krivitsky reference can indeed be a strong proxy for the latent dynamic process.  However, it should also be noted that the sensitivity to time scale separation is greater for $\bar{t}_\Delta$ than $\bar{d}$, and this sensitivity increases as $N$ grows and $P$ shrinks (i.e., as the number of foci grows).  Large systems with a slower underlying migration rate will tend to show excess triangulation that will not be captured by the reference model, and that will have to be modeled directly.

\section{Discussion} \label{sec_discuss}

Our examination of the conditions in which the Krivitsky reference model is expected to emerge from an underlying dynamic process raises a number of additional observations.  Here, we briefly comment on two further issues: the implications of this work for the importance of studying time scales in social phenomena, and the generalization of these results to latent spatial dynamics.

\subsection{The Importance of Time Scales}

While the importance of identifying the characteristic time scales on which different processes unfold is well-known in the physical and biological sciences, it has been less stressed within the social sciences.  As \citet{butts:s:2009} has observed, the relative time scales of edge turnover and network-based processes such as diffusion can give rise to radically different outcomes even where (time-aggregated) topology is preserved.  Relatedly, \cite{lerner.et.al:jmp:2013} have shown that the extent to which discretized network dynamics are affected by simultaneous dependence among edges is critically dependent upon the pace of network change relative to the time scale of measurement, with important implications for what must be considered when modeling network evolution.  Here, we have demonstrated another aspect of this phenomenon: social dynamics occurring on a time scale much faster than that of network formation can still leave their mark on the the networks they influence, despite being ephemeral from the point of view of any given relationship.  Even where such dynamics cannot be directly observed, we may in some cases be able to model their effects, and they may indeed provide a useful tool for intepreting structural biases (e.g., constant mean degree scaling) that might otherwise prove obscure. 

At the same time, the fact that ``fast'' processes can leave such systematic fingerprints on more slowly evolving degrees of freedom underscores the importance of identifying such processes, and characterizing them to the extent feasible.  That there is at this time no systematic reference or resource for assessing even the typical time scales of common social processes (much less their typical behavior) is in this regard a significant concern.  Developing a set of standardized time scales (and, where possible, reference models) for common social processes would seem to be an important aim for future research.

\subsection{Generalization to Spatial Processes}

While we have focused here on abstract foci that could have organizational, spatial, or other interpretations, further theoretical leverage can be obtained by endowing the foci with additional structure.  The most obvious example would be the case of latent spatial processes, in which individuals are assumed to migrate in an (unobserved) space and to be able to form ties only when sufficiently close to others.  Although a more elaborate treatment of this problem is possible, the following sketches one simple approach.

Let us consider a dynamic process in which individual nodes are confined to some region of volume $V$, defining $\rho = N/V$ to be the population density.\footnote{We make no particular assumptions here about the dimension of the space, or the metric with respect to which the volume is defined.}  We approximate the spatial constraint on tie formation by partitioning the region into $M$ uniform voxels of volume $v=V/M$, and imposing the condition that ties may only be formed by vertices occupying the same voxel.  As previously, we assume that vertices migrate at random between voxels at some rate $r_m$, with the migration process being such that (1) the equilibrium distribution of vertices in voxels is uniform and (2) the migration rate is sufficiently high that voxel membership is randomized on the time scale of the edge formation and dissolution processes.  In that case, we may invoke the result of Equation~\ref{e_psi} to obtain 
\begin{align}
\psi &= \log \frac{r_f v \rho}{r_\ell} - \log (V \rho), \nonumber \\
     &= \log \frac{r_f v }{r_\ell} - \log V. \label{e_psi_spa}
\end{align}
Thus, our corresponding ERGM family has reference measure $h(g) = V^{-t_e(g)}$ and edge parameter $\theta_e = \log \tfrac{r_f v }{r_\ell}$.  Interestingly, population density does not appear in the ERGM reference model: this is because it affects the system only through the number of vertices, and not through their geometry.  (See e.g. parallel results for the cross-sectional spatial models of \citet[pp84--86]{butts.et.al:sn:2012}.)  Increasing the population density under this model linearly increases $\bar{d}$ without changing the reference measure or parameter interpretation, just as increasing $N$ does in the conventional (non-sparse) ERGM family.  By contrast, increasing the system volume will increase $N$ \emph{without} increasing the mean degree, so long as the interaction volume ($v$) remains constant.

What if the interaction volume \emph{were} to change?  Where we wish to treat $v$ as exogenously varying, is useful to absorb it into the reference, i.e., moving it from the left-hand term of Equation~\ref{e_psi_spa} to the right-hand term.  This in turn gives us reference measure $h(g) = (v/V)^{t_e(g)}$ and changes the interpretation of $\theta_e$ to $r_f/r_\ell$.  Where $V$ is held constant, it can in turn be folded back into $\theta_e$, giving the constant system volume/varying interaction volume reference $h(g) = v^{t_e(g)}$ with corresponding parameter $\theta_e=r_f/(V r_\ell)$.  Although such changes are purely aesthetic when all elements of the model are known, they are important in inferential settings where $\theta_e$ is estimated from data.  Here, the choice of reference alters what is estimated, its units, and its interepretation, and generalization to new cases must take this into account.  Decompositions like those of Equation~\ref{e_psi_spa} can help to clarify the nature of such parameterization decisions.

Finally, we note that the above can be used to derive scaling relations when a suitable choice of geometry is made.  For instance, let us assume that our system occupies a hypercube of linear dimension $L$ in a $k$-dimensional Euclidean space, and that our voxels are defined by similar hypercubes of linear dimension $l$.  Now we may rewrite Equation~\ref{e_psi_spa} as
\[
\psi = \log \frac{r_f l^k }{r_\ell} - \log L^k,
\]
with $\theta_e$ and $h$ defined accordingly.  From this we can obtain scaling relations for the linear dimension of the interaction voxel or the system as a whole, or even for the dimensionality of the system.  One can envison this being particularly useful when applying models to relations among non-human animals, or to physical systems.

\section{Conclusion} \label{sec_conclusion}

We have here demonstrated that a common reference family for sparse graph models can be obtained as the result of an unobserved dynamic process in which individuals can form ties only when co-located in an abstract space of foci.  The ERGM reference model is a good proxy when the movement of vertices through foci is effectively time-scale separated from the network evolution process; simulation experiments suggest that a migration rate of one to two orders of magnitude faster than the edge formation rate is sufficient for this approximation to hold.  This is consistent with an interpretation of foci as transiently occupied settings within which vertices may encounter each other (e.g., gathering places, organizations, etc.), with such encounters being necessary for ties to form.  As we have shown, this framework can be straightforwardly generalized to vertices undergoing random motion in a region of finite volume, with the constraint that ties can only be formed when vertices are proximate.  This interpretation may be particularly useful for models of non-human animals or other systems (e.g., interactions among mobile devices).

Although our development provides one process-oriented interpretation of the Krivitsky reference model, it should be stressed that this is not the only interpretation possible.  Moreover, the observation of so-called ``power law densification'' \citep{leskovec.et.al:kdd:2007} in some systems suggests that the Krivitsky reference is not always ideal; the conditions necessary to recapitulate more general reference models like those of \citet{butts.almquist:jms:2015} are currently unknown.  Elaborating the relationships between dynamic processes and their associated reference families would thus seem to be a ripe area for further research.

\bibliography{ctb}


\end{document}